\documentclass[final,3p,times]{elsarticle}

\usepackage{graphicx,amsmath,amssymb,txfonts}
\usepackage{algorithm}
\usepackage{algorithmic}
\usepackage{subfigure,color}
\usepackage[colorlinks]{hyperref}
\usepackage{hypernat}
\usepackage{multirow}

\graphicspath{{Figs/}}

\biboptions{sort&compress}

\makeatletter

\newcommand{\Rmnum}[1]{\expandafter\@slowromancap\romannumeral #1@}
\makeatother

\journal{Applied Mathematics and Computation}

\newtheorem{theorem}{\textbf{Theorem}}

\newtheorem{lemma}{\textbf{Lemma}}

\begin{document}

\begin{frontmatter}



\title{A linearly implicit conservative  difference scheme for the generalized Rosenau-Kawahara-RLW equation}


\author[Tongji]{Dongdong He\corref{cor}}

\author[CSU,UNCC]{Kejia Pan}

\cortext[cor]{Corresponding author. E-mail address:dongdonghe@tongji.edu.cn(Dongdong He)}
\address[Tongji]{School of Aerospace Engineering and Applied Mechanics,
Tongji University, Shanghai 200092, China.}
\address[CSU]{School of Mathematics and Statistics, Central South University, Changsha 410083, China.}
 \address[UNCC]{Department of Mathematics and Statistics, University of North Carolina at Charlotte, Charlotte, NC 28223, USA}

\begin{abstract}
This paper concerns the numerical study for the generalized Rosenau-Kawahara-RLW equation obtained by coupling the generalized Rosenau-RLW equation and the generalized Rosenau-Kawahara equation.
We first derive the energy conservation law of the equation, and then develop a three-level linearly implicit difference scheme for solving the equation.
We prove that the proposed scheme is energy-conserved, unconditionally stable and  \textcolor{blue}{second-order accurate both in
time and space variables}. Finally, numerical experiments are carried out to confirm the energy
conservation, the convergence rates of the scheme and effectiveness for long-time simulation.
\end{abstract}


\end{frontmatter}
\textcolor{blue}{{\bf AMS subject classifications:}  35Q53, 65M06}

\section{Introduction}

The nonlinear wave is one of the most important scientific research areas. During the past several decades, many scientists developed different mathematical models to explain the wave behavior, such as the KdV equation, the RLW equation, the Rosenau equation, and many others. In the following, we give a short review of these important wave models.

The well known KdV equation
\begin{equation}\label{KdV}
u_{t}+u_{xxx}+6uu_x=0,
\end{equation}
was first introduced by  Boussinesq~\cite{Boussinesq} in 1877 and  rediscovered by Diederik Korteweg and Gustav de Vries~\cite{Korteweg} in 1895. Since then, there are a lot of studies  on this equation and its variational form. Here we just mention some of the recent work.  Kudryashov~\cite{Kudryashov} reviewed the travelling wave solutions for the KdV and the
KdV-Burgers equations proposed by~\cite{Wazzan}, Biswas~\cite{BiswasKdV} studied the solitary wave solution for KdV equation with power-law nonlinearity and time-dependent coefficients,  while Wang et al.~\cite{WangGW} investigated the solitons, shock waves for the potential KdV equation. In addition to the theoretical studies, readers can refer to~\cite{Dehghan,Vaneeva} for the numerical simulations of the KdV equation and the generalized KdV equation.

The regularized long-wave (RLW) equation (also known as  Benjamin-Bona-Mahony equation)
\begin{equation}\label{GRLW}
u_{t}+u_{x}+uu_{x}-u_{xxt}=0,
\end{equation}
was first proposed as a model for small-amplitude long wave of water in a channel by Peregrine~\cite{Peregrine1,Peregrine2}.  \textcolor{blue}{The regularized long-wave (RLW) equation and its different variational forms were well studied  both theoretically  and numerically in the literature. Readers can refer to ~\cite{Chegini, ZhangRLW, Roshan,Dehghan1,Dehghan2,Dehghan3,Dehghan4} for numerical studies  and ~\cite{BiswasRLW, Mohebbi, Song} for theoretical studies.}

When study the compact discrete systems, the well-known KdV equation can not describe the wave-wave
and wave-wall interactions. To overcome the shortcoming of KdV
equation, Rosenau proposed the following so-called Rosenau
equation~\cite{Rosenau, Rosenau2}:
\begin{equation}\label{Rosenaueq}
u_{t}+u_{xxxxt}+u_{x}+uu_x=0.
\end{equation}
The existence and uniqueness of the solution for the Rosenau equation were theoretically proved by~\cite{Park}. Besides the theoretical analysis,  numerical studies of equation (\ref{Rosenaueq}) also exist in the literature, see~\cite{Chung1, Chung2, Omrani, Manickam, Choo} and  references
therein.

For further consideration of nonlinear  waves, the term $-u_{xxt}$ is included in the Rosenau equation. The resulting  equation is usually called the Rosenau-RLW equation~\cite{ Pan1, Pan3, Atouani}:
\begin{equation}
u_{t}-u_{xxt}+u_{xxxxt}+u_{x}+uu_x=0.\label{RLWR}
\end{equation}
The above equation was further extended  into the generalized Rosenau-RLW equation~\cite{Pan2, Zuo2, Mittal}:
\begin{equation}
u_{t}-u_{xxt}+u_{xxxxt}+u_{x}+u^mu_x=0,\label{GRLWR}
\end{equation}
where $m\geq 1$ is a positive integer. The initial boundary value problems for the Rosenau-RLW equation and the  generalized Rosenau-RLW equation have been well studied numerically in the past
years~\cite{Pan1, Pan3, Atouani,  Pan2, Zuo2, Mittal}. For example, Pan et al.~\cite{Pan1} developed a three-level linearly implicit conservative scheme for the usual Rosenau-RLW equation~(\ref{RLWR}), where the method~\textcolor{blue}{achieves second-order accuracy both in time and space variables}. Zuo et al.~\cite{Zuo2} proposed
a \textcolor{blue}{second-order}  Crank-Nicolson finite difference scheme for the generalized Rosenau-RLW equation~(\ref{GRLWR}), however, the method is a nonlinear scheme,  iterations are needed at each time step.  Later on, Pan et al.~\cite{Pan2} used a three-level linearly implicit conservative scheme for the same generalized Rosenau-RLW equation~(\ref{GRLWR}), where the method is a linear scheme  and is also \textcolor{blue}{second-order} convergent  both in time and space variables.

On the other hand, to consider
another  behavior  of nonlinear waves, the viscous term $u_{xxx}$ needs to be included in the Rosenau equation (\ref{Rosenaueq}). The resulting  equation is usually called the
Rosenau-KdV equation~\cite{Esfahani, Saha, Razborova0, Zuo1, Hu1}:
\begin{equation}
u_{t}+u_{xxxxt}+u_{xxx}+u_{x}+uu_x=0.\label{RKdV}
\end{equation}
For numerical investigations, Hu et al.~\cite{Hu1} proposed a \textcolor{blue}{second-order} linear conservative finite difference method for the Rosenau-KdV equation. However, numerical methods for the initial-boundary value problem of the Rosenau-KdV equation have not been studied widely.

\textcolor{blue}{By coupling the above Rosenau-RLW equation and Rosenau-KdV equation, one can obtain the following Rosenau-KdV-RLW equation~\cite{Wongsaijai,Razborova1,Razborova2,Razborova3,Sanchez2015},
\begin{equation}
u_{t}-u_{xxt}+u_{xxxxt}+u_{xxx}+u_{x}+uu_x=0. \label{RKR}
\end{equation}
For numerical investigation, Wongsaijai et al.~\cite{Wongsaijai} proposed a three level implicit conservative  finite difference method for the above Rosenau-KdV-RLW equation. For theoretical  studies,  the solitary wave, shock waves, conservation laws as well as the asymptotic  behavior for the  Rosenau-KdV-RLW equation with  power law nonlinearity were studied by~\cite{Razborova1,Razborova2,Razborova3,Sanchez2015}, where the power law nonlinearity means the the last term in the left-hand side of equation (\ref{RKR}) is replaced by a general nonlinear term $(u^p)_x$ and $p$ is any positive integer.}

\textcolor{blue}{In addition, the following   Kawahara equation
\begin{equation}\label{GKE}
u_{t}+u_{x}+uu_{x}+ u_{xxx}-u_{xxxxx}=0,
\end{equation}
arose  in the theory of shallow water waves with surface tension~\cite{Kawahara}. Equation (\ref{GKE}) is called the modified  Kawahara equation if the third nonlinear term in the left-hand side is replaced by $u^2u_x$. There is a wide range of literature on the numerical investigations and theoretical studies for the usual Kawahara equation and the modified Kawahara equation. For theoretical aspects, some periodic and solitary wave solutions for both the Kawahara equation and the modified Kawahara equation are provided in~\cite{Sirendaoreji, Wazwaz, Yusufoglu}. In addition to the theoretical studies, readers can refer to~\cite{Polat, Jin, Dehghan5} for the numerical studies of the Kawahara equation and the modified  Kawahara equation.}

As one more step consideration of the nonlinear wave, Zuo~\cite{Zuo1}
obtained the Rosenau-Kawahara equation by adding another viscous
term  $-u_{xxxxx}$ to the Rosenau-KdV equation (\ref{RKdV}),  and studied
the solitary solution and periodic solution of the Rosenau-Kawahara equation. The  Rosenau-Kawahara equation is given as follows~\cite{Zuo1}:
\begin{equation}
u_{t}+u_{x}+uu_x+u_{xxx}+u_{xxxxt}-u_{xxxxx}=0.\label{RK}
\end{equation}
For numerical study, Hu et al.~\cite{Hu2} proposed a two level nonlinear Crank-Nicolson scheme and another three-level implicit linear  conservative finite difference scheme for the Rosenau-Kawahara equation, both methods are proved to be \textcolor{blue}{second-order} convergent. And Biswas~\cite{Biswas} investigated the solitary solution and the two invariance of following equation,
\begin{equation}
u_{t}+au_{x}+bu^{m}u_{x}+cu_{xxx}+\lambda u_{xxxxt}-\nu u_{xxxxx}=0,\label{GRK}
\end{equation}
where $a, b, c, \nu, \alpha,\lambda$ are real constants, $m$ is a positive integer, which indicates the power
law nonlinearity. Equation (\ref{GRK}) is referred as the generalized  Rosenau-Kawahara equation~\cite{Hu2}. As far as we know, besides the analytical solitary solution obtained by~\cite{Biswas}, there is no numerical study for this generalized  Rosenau-Kawahara equation.

\textcolor{blue}{Besides all the literature mentioned in the above, there are also a lot of related references which used both analytical and numerical methods to study  different nonlinear models. For example, Triki et al.~\cite{Triki2014} discussed the solitary solution for the Gear-Grimshaw model, Bhrawy et al.~\cite{Bhrawy2014} provided solitons, cnoidal waves and snoisal waves for the Whitham-Broer-Kaup system, while Shokri et al.~\cite{Dehghan6} numerically solved the two-dimensional complex Ginzburg-Landau equation by using the meshless radial basis functions method, and Dehghan et al.~\cite{Dehghan7} applied a semi-analytical method for solving the Rosenau-Hyman equation. Readers can refer to~\cite{ Triki2014,Bhrawy2014,Ebadi2013,Bhrawy20133} for more discussions on analytically finding solitons and other wave solutions for different  nonlinear models, and refer to~\cite{Dehghan6,Dehghan7,Dehghan8} for more discussions on numerical studies.}

In this paper, we will study the following generalized Rosenau-Kawahara-RLW equation:
\begin{equation}\label{KS1}
u_{t}+au_{x}+bu^{m}u_{x}+cu_{xxx}-\alpha u_{xxt}+\lambda u_{xxxxt}-\nu u_{xxxxx}=0, \quad  x_{l}\leq x \leq x_{r},\quad 0\leq t\leq T.
\end{equation}
which is obtained by coupling the generalized Rosenau-RLW equation (\ref{GRLWR}) and the generalized Rosenau-Kawahara equation (\ref{GRK}). \textcolor{blue}{Analytical solitary solutions of the above equation (\ref{KS1}) were recently obtained by~\cite{Zuon} through the sech ansatz method.}

\textcolor{blue}{In this work, we shall consider the  numerical simulation of above equation (\ref{KS1}) with initial condition}
\begin{equation}\label{KS2}
u(0,x)=u_{0}(x),  \quad  x_{l} \leq x \leq x_{r},
\end{equation}
and boundary conditions
\begin{align}\label{KS3}
&u(x_l,t)=u(x_r,t)=0,\quad u_x(x_l,t)=u_x(x_r,t)=0,\quad u_{xx}(x_l,t)=u_{xx}(x_r,t)=0,   \quad 0\leq t\leq T,
\end{align}
where $m$ is a positive integer, $a, b, c, \nu$ are all real constants, while $\alpha,\lambda$ are set to be positive constants, $x_l$ is a large negative number  and  $x_r$ is a large positive number. In this paper, we only discuss the solitary  solution of equation (\ref{KS1}). By solitary wave assumptions, the solitary solution and its derivatives have the following asymptotic values: $u\rightarrow 0$ as $x\rightarrow\pm\infty$, and $\frac{\partial^n u}{\partial x^n}\rightarrow 0$ as  $x\rightarrow\pm\infty$, for $n\geq1$. Thus, the boundary conditions (\ref{KS3}) are meaningful for the solitary solution of equation (\ref{KS1}). And we assume that the wave peak  is initially located at $x=0$, and  $x_l, x_r$, which are large numbers,  are used to assure that the solitary wave peak always locates inside the domain $[x_l, x_r]$ during the time interval $[0,T]$. Similar set-ups are used in~\cite{Hu1,Hu2,Wongsaijai}.

Note that equation (\ref{KS1}) reduces to the generalized Rosenau-RLW equation  when $c=0, \nu=0$. Specially,  it reduces to the usual Rosenau-RLW equation (\ref{RLWR}) when $a=m=\alpha=\lambda=1, b=2, c=0, \nu=0$. And equation (\ref{KS1}) reduces to the generalized Rosenau-Kawahara equation when $\alpha=0$.  Specially,  it reduces to the usual Rosenau-Kawahara equation  when $a=m=c=\lambda=\nu=1, b=2, \alpha=0$.

For numerically solving differential equations, the total accuracy of a particular method is affected not only the order of accuracy of the method, but also other factors.  The conservative property of the method is another factor that has the same or possibly even more impact on results.  For example, one successful and active research is to construct structure-preserving schemes (or called symplectic
schemes) for the ODE systems (see~\cite{Hairer} and the references therein). Better solutions can be expected from numerical schemes which have effective conservative properties rather than the ones which have nonconservative properties~\cite{Ham}.  And Li et al.~\cite{Li} even pointed out that in some areas, the ability to preserve some invariant properties of the original differential equation is
a criterion to judge the success of a numerical simulation.

For the initial boundary value problem (\ref{KS1})-(\ref{KS3}), we will show that it satisfies a fundamental energy conservation law in the next section. In addition, the wave equation is highly nonlinear due to the third nonlinear term in the left-hand side of the equation~(\ref{KS1}). When considering the finite difference scheme for the equation (\ref{KS1}), the usual Crank-Nicolson scheme will lead to a nonlinear scheme with heavy computation while other standard linearized discretizations for the nonlinear term, e.g., one step Newton's method or a \textcolor{blue}{second-order} extrapolation method, will loss the the  energy conservative property.  An ideal scheme should have relative less computational cost, can preserve energy, be unconditionally stable and maintain \textcolor{blue}{second-order} accuracy.

In  this paper, a three-level linearly implicit  finite difference method for the initial value problem (\ref{KS1})-(\ref{KS3}) will be presented. The fundamental energy conservation is preserved by the presented numerical scheme. The existence and uniqueness of the numerical solution are also proved. Moreover, numerical analysis shows that the method is \textcolor{blue}{second-order  convergent  both in time and space variables}, and  the method is unconditionally stable. Numerical results confirm well with the theoretical results.

The rest of the paper is organized  as follows:
 Section~\ref{section2} shows the energy conservation of the initial boundary value problem (\ref{KS1})-(\ref{KS3}).
 Section~\ref{section3}  gives the detail description of the  three-level linearly implicit finite difference method, the proof for the discrete conservative property,  the existence and uniqueness  as well as  the convergence and stability of the numerical solution. Numerical results are shown in section~\ref{section4}. Conclusions are provided in  the final section.

\section{Conservative property}\label{section2}
\textcolor{blue}{Equations (\ref{KS1})-(\ref{KS3}) satisfy} the following energy conservative property.
 \begin{theorem}\label{theorem1}
  Suppose $u_0\in C^{7}_{0}[x_l,x_r]$, then the solution of (\ref{KS1})-(\ref{KS3}) satisfies the energy conservation law:
 \begin{align}\label{energy}
E(t)=\int^{x_r}_{x_l}u^2(x,t)+\alpha u^2_{x}(x,t)+\lambda u^2_{xx}(x,t)dx=\int^{x_r}_{x_l}u^2(x,0)+\alpha u^2_{x}(x,0)+\lambda u^2_{xx}(x,0)dx=E(0),
\end{align}
for any $t\in[0,T]$, where $C^{7}_{0}[x_l,x_r]$ is the set of functions which are seventh  order continuous differentiable in the interval $[x_l, x_r]$ and have compact supports inside $(x_l, x_r)$.
 \end{theorem}

 {\bf Proof.} Multiplying (\ref{KS1}) by $2u$ and integrating  over the interval $[x_l, x_r]$, one get
 \begin{align}\label{proofe}
\frac{d}{dt}\int^{x_r}_{x_l}u^2dx&+a\int^{x_r}_{x_l}2uu_xdx+b\int^{x_r}_{x_l}2u^{m+1}\left(u\right)_xdx+c\int^{x_r}_{x_l}2uu_{xxx}dx\nonumber\\
&-\alpha\int^{x_r}_{x_l}2uu_{xxt}dx+\lambda\int^{x_r}_{x_l}2uu_{xxxxt}dx-\nu\int^{x_r}_{x_l}2uu_{xxxxx}dx=0.
\end{align}
Using the integration by parts, one can easily obtain
 \begin{align}
&\int^{x_r}_{x_l}uu_xdx=\int^{x_r}_{x_l}udu=\frac{1}{2}\left(u^2(x_r,t)-u^2(x_l,t)\right)=0,\nonumber\\
&\int^{x_r}_{x_l}u^{m+1}\left(u\right)_xdx=\int^{x_r}_{x_l}u^{m+1}du=\frac{1}{m+2}u^{m+2}\big|^{x_r}_{x_l}=0,\nonumber\\
&\int^{x_r}_{x_l}uu_{xxx}dx=uu_{xx}\big|^{x_r}_{x_l}-\int^{x_r}_{x_l}u_{xx}du=-\int^{x_r}_{x_l}u_{xx}u_xdx=-\int^{x_r}_{x_l}u_xdu_x=-\frac{1}{2}u^2_x\big|^{x_r}_{x_l}=0,  \nonumber\\
&\int^{x_r}_{x_l}uu_{xxt}dx=uu_{xt}\big|^{x_r}_{x_l}-\frac{1}{2}\frac{d}{dt}\int^{x_r}_{x_l}u^2_{x}dx=-\frac{1}{2}\frac{d}{dt}\int^{x_r}_{x_l}u^2_{x}dx,\nonumber\\
&\int^{x_r}_{x_l}uu_{xxxxt}dx=\int^{x_r}_{x_l}udu_{xxxt}=uu_{xxxt}\big|^{x_r}_{x_l}-\int^{x_r}_{x_l}u_xu_{xxxt}dx=-u_xu_{xxt}\big|^{x_r}_{x_l}+\frac{1}{2}\frac{d}{dt}\int^{x_r}_{x_l}u^{2}_{xx}dx=\frac{1}{2}\frac{d}{dt}\int^{x_r}_{x_l}u^{2}_{xx}dx,\nonumber\\
&\int^{x_r}_{x_l}uu_{xxxxx}dx=uu_{xxxx}\big|^{x_r}_{x_l}-\int^{x_r}_{x_l}u_{xxxx}u_xdx=-u_xu_{xxx}\big|^{x_r}_{x_l}+\int^{x_r}_{x_l}u_{xxx}u_{xx}dx=\frac{1}{2}u^{2}_{xx}\big|^{x_r}_{x_l}=0,\nonumber
\end{align}
where the boundary conditions (\ref{KS3}) are used.

Thus, only the first, fifth and sixth term in the left hand side of (\ref{proofe}) are nonzero, all other terms are  zero. This yields,
 \begin{align}
\frac{d}{dt}\int^{x_r}_{x_l}u^2+\alpha u^2_{x}+\lambda u^2_{xx}dx=0.
\end{align}
Therefore,
 \begin{align}
E(t)=\int^{x_r}_{x_l}u^2(x,t)+\alpha u^2_{x}(x,t)+\lambda u^2_{xx}(x,t)dx=\int^{x_r}_{x_l}u^2(x,0)+\alpha u^2_{x}(x,0)+\lambda u^2_{xx}(x,0)dx=E(0),
\end{align}
for any $t\in [0,T]$. This completes the proof.

%

\section{Numerical method}\label{section3}

\subsection{Numerical scheme}
In this section, we give a complete description of our numerical method for the initial value problem (\ref{KS1})-(\ref{KS3}).  We first describe our solution domain and its grid. The computational domain is defined as $\Omega=\{(x,t) | x_l\leq x\leq x_r,\ 0\leq t\leq T\}$, which is covered by a uniform grid $\Omega_{h}=\{(x_i,t_n)|x_i=x_l+ih,\ t_n=n\tau,\ i=0, \cdots, M,\ n=0, \cdots, N\}$, with spacing $h=\frac{x_r-x_l}{M},\ \tau=\frac{T}{N}$. We denote $U^{n}_{i}$ as the numerical approximation of $u(x_i,t_n)$ and $Z^{0}_{h}=\{U=(U_i)|U_{-1}=U_{0}=U_{1}=U_{M-1}=U_{M}=U_{M+1}=0,\ i=-1,0,1, \cdots, M-1,M, M+1\}$. For convenience, the difference operators, inner product and norms  are defined as follows:
\begin{align*}
&\bar{U}^{n}_{i}=\frac{U^{n+1}_{i}+U^{n-1}_{i}}{2}, \quad (U^{n}_i)_{t}=\frac{U^{n+1}_i-U^{n-1}_i}{2\tau}, \quad (U^{n}_{i})_{x}=\frac{U^{n}_{i+1}-U^{n}_i}{h},\quad (U^{n}_{i})_{\bar{x}}=\frac{U^{n}_{i}-U^{n}_{i-1}}{h},\\
&(U^{n}_{i})_{\hat{x}}=\frac{U^{n}_{i+1}-U^{n}_{i-1}}{2h},\quad (U^n,V^n)=h\sum^{M-1}_{j=1}U^{n}_iV^{n}_i,\quad \parallel U^n\parallel^2=(U^n,U^n),\quad \parallel U^n\parallel_{\infty}=\max_{1\leq i \leq M-1}|U^{n}_i|.
\end{align*}
where $U^n=(U^n_{-1},\cdots,U^n_{M+1})$ is the numerical solution at time $t_n=n\tau$.

The essential of our  scheme is that the third term  in the left-hand side of (\ref{KS1}) is rewritten and discretized as
\begin{align}\label{equation1}
bu^mu_x=\frac{b}{m+2}(u^mu_x+(u^{m+1})_x)\approx\frac{b}{m+2}[(U^n_i)^{m}(\bar{U}^{n}_{i})_{\hat{x}}+((U^n_i)^{m}\bar{U}^{n}_{i})_{\hat{x}}],
\end{align}
and
which is  a \textcolor{blue}{second-order} approximation around $(x_i=x_l+ih,t_{n}=n\tau)$. And all other terms  of left-hand side of (\ref{KS1}) are discretizated  by using  the standard \textcolor{blue}{second-order} central difference method around $(x_i=x_l+ih,t_{n}=n\tau)$ .

The \textcolor{blue}{detailed}  numerical scheme is as follows:
 \begin{align}\label{scheme1}
 (U^n_i)_t+a(\bar{U}^{n}_{i})_{\hat{x}}&+\frac{b}{m+2}[(U^n_i)^{m}(\bar{U}^{n}_{i})_{\hat{x}}+((U^n_i)^{m}\bar{U}^{n}_{i})_{\hat{x}}]+c(\bar{U}^{n}_{i})_{x\bar{x}\hat{x}}-\alpha(U^{n}_i)_{x\bar{x}t}
 \nonumber \\&+\lambda(U^{n}_{i})_{xx\bar{x}\bar{x}t}-\nu(\bar{U}^{n}_{i})_{xx\bar{x}\bar{x}\hat{x}}=0,\ \  2\leq i\leq M-2,\ \ 1\leq n \leq N-1.
 \end{align}
and
  \begin{align}\label{scheme2}
 U^{0}_{i}=u_0(x_i), \  0\leq i\leq M,
  \end{align}
 \begin{align}\label{scheme3}
 U^{j}_{0}= U^{j}_{M}=0,  \ (U^{j}_{0})_{\hat{x}}= (U^{j}_{M})_{\hat{x}}=0,\ (U^{j}_{0})_{x\bar{x}}= (U^{j}_{M})_{x\bar{x}}=0,  \ 0\leq j \leq N.
 \end{align}
 Obviously, the above conditions (\ref{scheme3}) will give  $U^j_1=U^j_{M-1}=0$  on two internal points and $U^j_{-1}=U^j_{M+1}=0$ on two fictitious  points, for any $0\leq j\leq N$. Thus, $U^j \in Z^{0}_{h}$, for any $0\leq j\leq N$. In addition,  we can see that (\ref{scheme1}) is a three-level linearly implicit scheme and  the coefficient matrix of \textcolor{blue}{the linear system} (\ref{scheme1}) is banded. Therefore, the resulting \textcolor{blue}{linear system of equations}  can be solved efficiently using a direct linear solver, such as the LU decomposition method.

Since the scheme is a three-level method, to start the computation, we need to give the method for computation of $U^1$. The $U^1$ is computed through the following Crank-Nicolson scheme:
 \begin{align}\label{scheme4}
&\frac{U^1_i-U^0_i}{\tau}+a\left(\frac{U^1_i+U^0_i}{2}\right)_{\hat{x}}+\frac{b}{m+2}\left[\left(\frac{U^1_i+U^0_i}{2}\right)^{m}\left(\frac{U^1_i+U^0_i}{2}\right)_{\hat{x}}+\left(\left(\frac{U^1_i+U^0_i}{2}\right)^{m}\left(\frac{U^1_i+U^0_i}{2}\right)\right)_{\hat{x}}\right]\nonumber\\
&-\alpha\left(\frac{U^{1}_i-U^{0}_i}{\tau}\right)_{x\bar{x}}+c\left(\frac{U^1_i+U^0_i}{2}\right)_{x\bar{x}\hat{x}}+\lambda\left(\frac{U^1_i-U^0_i}{\tau}\right)_{xx\bar{x}\bar{x}}-\nu\left(\frac{U^1_i+U^0_i}{2}\right)_{xx\bar{x}\bar{x}\hat{x}}=0,
 \end{align}
 where it is a nonlinear  scheme and is \textcolor{blue}{second-order accurate both in time and space variables}.

 We point out that this discretization (\ref{equation1}) for the  third term  in the left-hand side of (\ref{KS1}) is specially designed so that the numerical scheme (\ref{scheme1})-(\ref{scheme3}) can guarantee the energy conservation (see the proof of theorem~\ref{theorem2} in the next subsection).  Readers can refer to references~\cite{Pan1,Hu2,Wongsaijai} for the similar kind of treatments for the nonlinear term $uu_x$ in their equations to achieve conservative schemes. In addition, our scheme (\ref{scheme1})-(\ref{scheme4}) is only nonlinear at first time step when computing $U^1$, and all other steps are linear. Thus, the computational cost is relatively cheap. In the following several sections of the paper,  we will show that our scheme is \textcolor{blue}{second-order accurate both in time and space variables}, can preserve the energy identity (\ref{energy}), and is unconditionally stable.

The following lemmas are well known results, which are essential for existence, uniqueness, convergence, and stability of our numerical solution. \textcolor{blue}{Lemma~\ref{Lemma1} and~\ref{Lemma2}} can be verified  through direct computation, lemma~\ref{Lemma3} is the discrete form of the Sobolev's embedding theorem, which can be found in lemma~\ref{Lemma1}, page 110  of \cite{Samarskii}.
 In the rest part of the paper, unless otherwise indicated, $C$ is the notation referring to a general positive constant, which may have different values in different \textcolor{blue}{contexts}.\\

\begin{lemma}\label{Lemma1}
  For any two mesh functions $U, V\in Z^{0}_{h}$, one  have
\begin{align}
&(U_{\hat{x}}, V)=-(U, V_{\hat{x}}), \quad (U_{x}, V)=-(U, V_{\bar{x}}), \quad (U_{x\bar{x}}, V)=-(U_x,V_x), \nonumber
\end{align}
 Furthermore, 
 \begin{align}
  (U, U_{xx\bar{x}\bar{x}})=\parallel U_{x\bar{x}}\parallel^2.\nonumber
 \end{align}
\end{lemma}

\begin{lemma}\label{Lemma2}
For any mesh function $U\in Z^{0}_{h}$, one have
 \begin{align}
  (U_{\hat{x}},U)=0,\quad (U_{x\bar{x}\hat{x}},U)=0,\quad (U_{xx\bar{x}\bar{x}\hat{x}},U)=0.\nonumber
 \end{align}
 \end{lemma}

\begin{lemma}\label{Lemma3}
(Discrete Sobolev's inequality (Lemma~\ref{Lemma1}, page 110  of~\cite{Samarskii}) For any mesh function $U\in Z^{0}_{h}$, one have
  \begin{align}
\parallel U\parallel_{\infty}\leq C\parallel U_{x}\parallel.\nonumber
  \end{align}
 \end{lemma}

\subsection{Discrete conservation}


 \begin{theorem}\label{theorem2}
 Suppose $u_0\in C^{7}_{0}[x_l,x_r]$, then the  solution of finite difference scheme (\ref{scheme1})-(\ref{scheme3}) \textcolor{blue}{satisfies} $\parallel U^{n}\parallel_{\infty}\leq C$, $\parallel U^{n}_{x}\parallel_{\infty}\leq C$, for any $0\leq n \leq N$, and also \textcolor{blue}{satisfies} the following discrete  energy conservation:
\begin{align}\label{disenergy}
 E^n&\triangleq\frac{\parallel U^{n+1}\parallel^2+\parallel U^{n}\parallel^2}{2}+\alpha\frac{\parallel U^{n+1}_{x}\parallel^2+\parallel U^{n}_{x}\parallel^2}{2}+\lambda\frac{\parallel U^{n+1}_{x\bar{x}}\parallel^2+\parallel U^{n}_{x\bar{x}}\parallel^2}{2}\nonumber\\
 &=\frac{\parallel U^{1}\parallel^2+\parallel U^{0}\parallel^2}{2}+\alpha\frac{\parallel U^{1}_{x}\parallel^2+\parallel U^{0}_{x}\parallel^2}{2}+\lambda\frac{\parallel U^{1}_{x\bar{x}}\parallel^2+\parallel U^{0}_{x\bar{x}}\parallel^2}{2}\triangleq E^0,
\end{align}
for any $0\leq n \leq N-1$, where $E^{n}$ is the discrete energy at time $t=(n+\frac{1}{2})\tau$.
 \end{theorem}


{\bf Proof}. Taking the conditions $U^j_{-1}=U^j_{0}=U^j_{1}=U^j_{M-1}=U^j_{M}=U^j_{M+1}=0$ ($0\leq j \leq N$) into account, and after computing the inner product of equation (\ref{scheme1}) with $\bar{U}^{n}$, i. e., $\frac{U^{n+1}+U^{n-1}}{2}$, we have
\begin{align}
&\frac{\parallel U^{n+1}\parallel^2-\parallel U^{n-1}\parallel^2}{4\tau}+a(\bar{U}^n_{\hat{x}},\bar{U}^n)+\frac{b}{m+2}\left((U^n)^{m}(\bar{U}^{n})_{\hat{x}}+((U^n)^{m}\bar{U}^{n})_{\hat{x}},\bar{U}^{n}\right)\\
&+c(\bar{U}^{n}_{x\bar{x}\hat{x}},\bar{U}^{n})-\alpha(U^n_{x\bar{x}t},\bar{U}^{n})+\lambda(U^n_{xx\bar{x}\bar{x}t},\bar{U}^{n})-\nu(\bar{U}^{n}_{xx\bar{x}\bar{x}\hat{x}},\bar{U}^{n})=0.
\end{align}
By using lemma~\ref{Lemma2}, we get
\begin{align}
(\bar{U}^{n}_{\hat{x}},\bar{U}^{n})=0, \quad (\bar{U}^{n}_{x\bar{x}\hat{x}},\bar{U}^{n})=0,\quad  (\bar{U}^{n}_{xx\bar{x}\bar{x}\hat{x}},\bar{U}^{n})=0.
\end{align}
Moreover,
  \begin{align}
&\left(\left[\left(U^n\right)^{m}\left(\bar{U}^{n}\right)_{\hat{x}}+\left(\left(U^n\right)^{m}\bar{U}^{n}\right)_{\hat{x}}\right],\bar{U}^{n}\right)\nonumber\\
&=\frac{1}{2}\sum^{M-1}_{i=1}\left[ (U^{n}_i)^{m}(\bar{U}^{n}_{i+1}-\bar{U}^{n}_{i-1})+(U^n_{i+1})^{m}\bar{U}^{n}_{i+1}-(U^n_{i-1})^{m}\bar{U}^{n}_{i-1}\right ]\bar{U}^{n}_i\nonumber\\
&=\frac{1}{2}\sum^{M-1}_{i=1}\left[ (U^{n}_i)^{m}\bar{U}^{n}_{i+1}\bar{U}^{n}_i+(U^n_{i+1})^{m}\bar{U}^{n}_{i+1}\bar{U}^{n}_i\right ]-\frac{1}{2}\sum^{M-1}_{i=1}\left[ (U^n_{i-1})^{m}\bar{U}^{n}_i\bar{U}^{n}_{i-1}+(U^{n}_i)^{m}\bar{U}^{n}_i\bar{U}^{n}_{i-1} \right]\nonumber\\
&=0.
 \end{align}
and
\begin{align}
(U^n_{xx\bar{x}\bar{x}t},\bar{U}^{n})=\frac{\parallel U^{n+1}_{x\bar{x}}\parallel^2-\parallel U^{n-1}_{x\bar{x}}\parallel^2}{4\tau}, \quad (U^n_{x\bar{x}t},\bar{U}^{n})=-\frac{\parallel U^{n+1}_{x}\parallel^2-\parallel U^{n-1}_{x}\parallel^2}{4\tau},
\end{align}
where boundary conditions (\ref{scheme2}) and (\ref{scheme3}) are used.

Thus,
\begin{align}
\parallel U^{n+1}\parallel^2-\parallel U^{n-1}\parallel^2+\alpha\parallel U^{n+1}_{x}\parallel^2-\alpha\parallel U^{n-1}_{x}\parallel^2+\lambda\parallel U^{n+1}_{x\bar{x}}\parallel^2-\lambda\parallel U^{n-1}_{x\bar{x}}\parallel^2=0,
\end{align}
for any $1\leq n\leq N-1$.  This is equivalent to

\begin{align}
&\frac{\parallel U^{n+1}\parallel^2+\parallel U^{n}\parallel^2}{2}+\alpha\frac{\parallel U^{n+1}_{x}\parallel^2+\parallel U^{n}_{x}\parallel^2}{2}+\lambda\frac{\parallel U^{n+1}_{x\bar{x}}\parallel^2+\parallel U^{n}_{x\bar{x}}\parallel^2}{2}=\nonumber\\
&\frac{\parallel U^{n}\parallel^2+\parallel U^{n-1}\parallel^2}{2}+\alpha\frac{\parallel U^{n}_{x}\parallel^2+\parallel U^{n-1}_{x}\parallel^2}{2}+\lambda\frac{\parallel U^{n}_{x\bar{x}}\parallel^2+\parallel U^{n-1}_{x\bar{x}}\parallel^2}{2},
\end{align}
for any $1\leq n\leq N-1$.
This further yields
\begin{align}\label{Energy_conservative}
E^n=E^0, \quad \textrm{for any}\ 1\leq n\leq N-1,
\end{align}
which is actually the energy conservation law (\ref{disenergy}).

Multiplying (\ref{scheme4}) both sides by $\frac{U^1_i+U^0_i}{2}$ and using the similar techniques as above, one can obtain
\begin{align}
\parallel U^{0}\parallel^2+\alpha\parallel U^{0}_{x}\parallel^2+\lambda\parallel U^{0}_{x\bar{x}}\parallel^2=\parallel U^{1}\parallel^2+\alpha\parallel U^{1}_{x}\parallel^2+\lambda\parallel U^{1}_{x\bar{x}}\parallel^2.
\end{align}

Thus, (\ref{Energy_conservative}) can be rewritten as
\begin{align}\label{Energy_conservative2}
E^n=\parallel U^{0}\parallel^2+\alpha\parallel U^{0}_{x}\parallel^2+\lambda\parallel U^{0}_{x\bar{x}}\parallel^2.
\end{align}

Since $u_0\in C^7_0[x_l, x_r]$ and the initial condition (\ref{scheme2}) are used in the numerical method,  the right-side of (\ref{Energy_conservative2}) is bounded. By assumptions, $\alpha, \lambda$ are positive constants,  therefore,
\begin{align}
 \parallel U^{n}_{x}\parallel\leq C,\ \   \parallel U^{n}_{x\bar{x}}\parallel\leq C,\quad  \textrm{for any}\ 0\leq n\leq N.
\end{align}

By using lemma~\ref{Lemma3}, we have $\parallel U^{n}\parallel_{\infty}\leq C$.

In addition, through direct computation, one can  verify that
\begin{align}
 \parallel U^{n}_{xx}\parallel= \parallel U^{n}_{x\bar{x}}\parallel, \quad  \textrm{for any}\ 0\leq n\leq N.
\end{align}

Thus,
\begin{align}
 \parallel U^{n}_{xx}\parallel\leq C, \quad  \textrm{for any}\ 0\leq n\leq N.
\end{align}

Again by using lemma~\ref{Lemma3}, we have $\parallel U^{n}_x\parallel_{\infty}\leq C$. This completes the proof.

\subsection{Existence and uniqueness}
 \begin{theorem}\label{theorem3}
 The finite difference scheme (\ref{scheme1})-(\ref{scheme3}) has a unique solution.
 \end{theorem}

\textbf{Proof}. To prove the \textcolor{blue}{theorem}, we proceed by the mathematical induction. Suppose $U^1,\cdots,U^{n} (1\leq n\leq N-1)$ are solved uniquely, we now consider the equation (\ref{scheme1}) for $U^{n+1}$.
Assume that  $U^{n+1,1},U^{n+1,2}$ are two solutions of (\ref{scheme1}) and let $W^{n+1}=U^{n+1,1}-U^{n+1,2}$,  then it is easy to verify that  $W^{n+1}$ satisfies the following equation:
 \begin{align}\label{unique}
 &\frac{1}{2\tau}W^{n+1}_i+\frac{a}{2}(W^{n+1}_i)_{\hat{x}}+\frac{b}{2(m+2)}\left[(U^n_i)^{m}(W^{n+1}_{i})_{\hat{x}}+((U^n_i)^{m}W^{n+1}_{i})_{\hat{x}}\right]\nonumber\\
 &+\frac{c}{2}(W^{n+1}_{i})_{x\bar{x}\hat{x}}-\frac{\alpha}{2\tau}(W^{n+1}_{i})_{x\bar{x}}+\frac{\lambda}{2\tau}(W^{n+1}_{i})_{xx\bar{x}\bar{x}}-\frac{\nu}{2}(W^{n+1}_{i})_{xx\bar{x}\bar{x}\hat{x}}=0.
 \end{align}

 Taking the inner product of (\ref{unique}) with $W^{n+1}$, we have
  \begin{align}\label{unique2}
 \frac{1}{2\tau}\parallel W^{n+1}\parallel^2+\frac{\alpha}{2\tau}\parallel W^{n+1}_{x}\parallel^2+ \frac{\lambda}{2\tau}\parallel W^{n+1}_{x\bar{x}}\parallel^2=0,
 \end{align}
 where
 \begin{align}\label{unique3}
&(W^{n+1}_{x\bar{x}},\quad W^{n+1})=-(W^{n+1}_{x},W^{n+1}_{x}),\quad (W^{n+1}_{xx\bar{x}\bar{x}},W^{n+1})=\parallel W^{n+1}_{x\bar{x}}\parallel^2,\quad (W^{n+1}_{xx\bar{x}\bar{x}\hat{x}},W^{n+1})=0,\nonumber\\ & (W^{n+1}_{x\bar{x}\hat{x}},W^{n+1})=0,\quad (W^{n+1}_{\hat{x}},W^{n+1})=0,\quad \left(\left[\left(U^n\right)^{m}\left(W^{n+1}\right)_{\hat{x}}+\left(\left(U^n\right)^{m}W^{n+1}\right)_{\hat{x}}\right],W^{n+1}\right)=0,
 \end{align}
 are used. The first five identities of (\ref{unique3}) are directly from lemma~\ref{Lemma1} and~\ref{Lemma2}, and the last one can be obtained  as follows:
  \begin{align}
&\left(\left[\left(U^n\right)^{m}\left(W^{n+1}\right)_{\hat{x}}+\left(\left(U^n\right)^{m}W^{n+1}\right)_{\hat{x}}\right],W^{n+1}\right)\nonumber\\
&=\frac{1}{2}\sum^{M-1}_{i=1}\left[ (U^{n}_i)^{m}(W^{n+1}_{i+1}-W^{n+1}_{i-1})+(U^n_{i+1})^{m}W^{n+1}_{i+1}-(U^n_{i-1})^{m}W^{n+1}_{i-1}\right ]W^{n+1}_i\nonumber\\
&=\frac{1}{2}\sum^{M-1}_{i=1}\left[ (U^{n}_i)^{m}W^{n+1}_{i+1}W^{n+1}_i+(U^n_{i+1})^{m}W^{n+1}_{i+1}W^{n+1}_i\right ]-\frac{1}{2}\sum^{M-1}_{i=1}\left[ (U^{n}_i)^{m}W^{n+1}_{i-1}W^{n+1}_i+(U^n_{i-1})^{m}W^{n+1}_{i-1}W^{n+1}_i\right ]\nonumber\\
&=0.
 \end{align}

 From (\ref{unique2}) and the definition of the $\parallel\cdot\parallel$-norm, one can see that (\ref{unique2}) has only a trivial solution. Thus, (\ref{scheme1}) determines $U^{n+1}$ uniquely.  This completes the proof.

\subsection{Convergence and stability}

Let $u(x,t)$ be the solution of problem (\ref{KS1})-(\ref{KS3}),  $U^n_i$ be the solution of the numerical schemes (\ref{scheme1})-(\ref{scheme3}), and  $u^n_i=u(x_i,t_n),\ e^n_i=u^n_i-U^n_i$, then the truncation error of the scheme (\ref{scheme1})-(\ref{scheme3}) can be obtained as follows:
\begin{align}\label{trucation}
r^n_i=&\left(e^n_i\right)_t+a(\bar{e}^{n}_i)_{\hat{x}}+\frac{b}{m+2}\left[(u^n_i)^{m}(\bar{u}^{n}_{i})_{\hat{x}}+\left(\left(u^n_i\right)^{m}\bar{u}^{n}_{i}\right)_{\hat{x}}-\left(U^n_i\right)^{m}\left(\bar{U}^{n}_{i}\right)_{\hat{x}}-\left(\left(U^n_i\right)^{m}\bar{U}^{n}_{i}\right)_{\hat{x}}\right]\nonumber\\
&+c\left(\bar{e}^{n}_{i}\right)_{x\bar{x}\hat{x}}-\alpha(e^{n}_{i})_{x\bar{x}t}+\lambda(e^n_i)_{xx\bar{x}\bar{x}t}-\nu\left(\bar{e}^{n}_{i}\right)_{xx\bar{x}\bar{x}\hat{x}},
\end{align}
where $\bar{e}^{n}=\frac{e^{n+1}+e^{n-1}}{2}$, $2\leq i\leq M-2$ and $1\leq n \leq N-1$.

Since all terms in (\ref{scheme1}) are the \textcolor{blue}{second-order} approximations of the corresponding terms in left-hand side of (\ref{KS1}) around $(x_i=x_l+ih, t_{n}=n\tau)$, by Taylor expansion, it can be easily obtained that $r^n_i=O(\tau^2+h^2)$ if $h, \tau\rightarrow 0$ and $u(x,t)\in C^{7,3}$,  where $C^{7,3}$ is the set of functions which are  seventh order continuous  differentiable in space and third  order continuous differentiable in time. \textcolor{blue}{This following lemma is a well known result}.\\

\begin{lemma}\label{Lemma4}
 (Discrete Gronwall's inequality).  Suppose that $w(k)$ and $\rho(k)$ are nonnegative functions while $\rho(k)$ is a non-decreasing function. If
\begin{align*}
w(k)\leq \rho(k)+C\tau \sum^{k-1}_{l=0}w(l), \forall\ k,
\end{align*}
then
\begin{align*}
w(k)\leq \rho(k)e^{C\tau k}, \forall\ k.
\end{align*}
\end{lemma}

\begin{theorem}\label{theorem4}
 Suppose $u_0\in C^{7}_{0}[x_l,x_r]$, and $u(x,t)\in C^{7,3}$, then the numerical solution $U^n$ of the finite difference scheme (\ref{scheme1})-(\ref{scheme3})  converges to the solution of the problem (\ref{KS1})-(\ref{KS3}) in the sense of  $\parallel\cdot\parallel_{\infty}$,  and the convergence rate is $O(\tau^2+h^2)$, i.e.,
 \begin{align}
\parallel u^n-U^n\parallel_{\infty}\leq C(\tau^2+h^2), \quad \textrm{for any}\ 2\leq n \leq N.
\end{align}
\end{theorem}

\textbf{Proof}.  Taking the inner product of (\ref{trucation}) with $2\bar{e}^{n}$, we have
\begin{align}\label{estimate0}
&\parallel e^{n+1}\parallel^2-\parallel e^{n-1}\parallel^2+\alpha(\parallel e^{n+1}_{x}\parallel^2-\parallel e^{n-1}_{x}\parallel^2)
+\lambda(\parallel e^{n+1}_{x\bar{x}}\parallel^2-\parallel e^{n-1}_{x\bar{x}}\parallel^2)\nonumber\\
&=2\tau\left[\left(r^n,2\bar{e}^{n}\right)-a\left(\bar{e}^{n}_{\hat{x}},2\bar{e}^{n}\right)-c\left(\left(\bar{e}^{n}\right)_{x\bar{x}\hat{x}},2\bar{e}^{n}\right)+\nu\left(\left(\bar{e}^{n}\right)_{xx\bar{x}\bar{x}\hat{x}},2\bar{e}^{n}\right)-\left(Q,2\bar{e}^{n}\right)-\left(R,2\bar{e}^{n}\right)\right],
\end{align}
where
\begin{align}
Q=\frac{b}{m+2}\left[\left(u^n\right)^{m}\left(\bar{u}^n\right)_{\hat{x}}-\left(U^n\right)^{m}\left(\bar{U}^n\right)_{\hat{x}}\right], \ R=\frac{b}{m+2}\left[\left(\left(u^n\right)^{m}\bar{u}^n\right)_{\hat{x}}-\left(\left(U^n\right)^{m}\bar{U}^n\right)_{\hat{x}}\right].
\end{align}

By using lemma~\ref{Lemma2}, we obtain
\begin{align}
\left(\bar{e}^{n}_{\hat{x}},2\bar{e}^{n}\right)=0, \quad \left((\bar{e}^{n})_{x\bar{x}\hat{x}},2\bar{e}^{n}\right)=0, \quad  \left((\bar{e}^{n})_{xx\bar{x}\bar{x}\hat{x}},2\bar{e}^{n}\right)=0.
\end{align}

\textcolor{blue}{From the notations introduced at the beginning of section~\ref{section3},  for any $0\leq n\leq N$, we have the following inequality}
\begin{align*}
\textcolor{blue}{\parallel e^{n}_{\hat{x}}\parallel^2=\sum^{M-1}_{j=1}\left(\frac{e_{j+1}-e_{j-1}}{2h}\right)^2h=\frac{1}{4}\sum^{M-1}_{j=1}\left(\frac{e_{j+1}-e_{j}}{h}+\frac{e_{j}-e_{j-1}}{h}\right)^2h\leq\frac{1}{2}\sum^{M-1}_{j=1}\left[\left(\frac{e_{j+1}-e_{j}}{h}\right)^2+\left(\frac{e_{j}-e_{j-1}}{h}\right)^2\right]h=\parallel e^n_{x}\parallel^2.}
\end{align*}
\textcolor{blue}{Thus,  for any $0\leq n\leq N$, we have }
\begin{align}\label{normin}
\textcolor{blue}{\parallel e^n_{\hat{x}}\parallel\leq \parallel e^n_{x}\parallel.}
\end{align}

In addition, we have
 \begin{align}\label{estimate1}
|(Q,2\bar{e}^{n})|&=\left|\frac{2b h}{m+2}\sum^{M-1}_{j=1}\left[(u^n_j)^{m}(\bar{u}^n_j)_{\hat{x}}-(U^n_j)^{m}(\bar{U}^n_j)_{\hat{x}}\right]\bar{e}^{n}_{j}\right|\nonumber\\
&=\left|\frac{2 b h}{m+b}\sum^{M-1}_{j=1}\left[(u^n_j)^{m}(\bar{u}^n_j)_{\hat{x}}-(u^n_j)^{m}(\bar{U}^n_j)_{\hat{x}}+(u^n_j)^{m}(\bar{U}^n_j)_{\hat{x}}-(U^n_j)^{m}(\bar{U}^n_j)_{\hat{x}}\right]\bar{e}^{n}_{j}\right|\nonumber\\
&=\left|\frac{2b h}{m+2}\sum^{M-1}_{j=1}(u^n_j)^{m}(\bar{e}^n_j)_{\hat{x}}\bar{e}^{n}_{j}+\frac{2 bh}{m+2}\sum^{M-1}_{j=1}\left[\left(u^n_j\right)^{m}-(U^n_j)^{m}\right](\bar{U}^n_j)_{\hat{x}}\bar{e}^{n}_{j}\right|\nonumber\\
&=\left|\frac{2b h}{m+2}\sum^{M-1}_{j=1}(u^n_j)^{m}(\bar{e}^n_j)_{\hat{x}}\bar{e}^{n}_{j}+\frac{2 bh}{m+2}\sum^{M-1}_{j=1}\left[\sum^{m-1}_{k=0}(u^n_j)^{m-1-k}(U^n_j)^{k}\right]e^n_j(\bar{U}^n_j)_{\hat{x}}\bar{e}^{n}_{j}\right|\nonumber\\
&\leq Ch\sum^{M-1}_{j=1}\left[|(\bar{e}^n_j)_{\hat{x}}|+|e^n_j|\right]|\bar{e}^{n}_{j}|\nonumber\\
&\leq C(\parallel e^{n+1}_{\hat{x}}\parallel^2+\parallel e^{n-1}_{\hat{x}}\parallel^2+\parallel e^{n+1}\parallel^2+\parallel e^{n}\parallel^2+\parallel e^{n-1}\parallel^2),\nonumber\\
&\leq C(\parallel e^{n+1}_x\parallel^2+\parallel e^{n-1}_x\parallel^2+\parallel e^{n+1}\parallel^2+\parallel e^{n}\parallel^2+\parallel e^{n-1}\parallel^2),
\end{align}
where theorem~\ref{theorem2}, Cauchy-Schwarz inequality  and inequality (\ref{normin}) are used.

Similarly, we have
 \begin{align}\label{estimate2}
|(R,2\bar{e}^{n})|&=\left|\frac{2bh}{m+2}\sum^{M-1}_{j=1}\left[\left((u^n_j)^{m}\bar{u}^n_j\right)_{\hat{x}}-\left((U^n_j)^{m}\bar{U}^n_j\right)_{\hat{x}}\right]\bar{e}^{n}_{j}\right|\nonumber\\
&=\left|\frac{2bh}{m+2}\sum^{M-1}_{j=1}\left[\left((u^n_j)^{m}\bar{u}^n_j\right)_{\hat{x}}-\left((u^n_j)^{m}\bar{U}^n_j\right)_{\hat{x}}+\left((u^n_j)^{m}\bar{U}^n_j\right)_{\hat{x}}-\left((U^n_j)^{m}\bar{U}^n_j\right)_{\hat{x}}\right]\bar{e}^{n}_{j}\right|\nonumber\\
&=\left|\frac{2bh}{m+2}\sum^{M-1}_{j=1}\left[\left((u^n_j)^{m}\bar{e}^n_j\right)_{\hat{x}}+\left(\sum^{m-1}_{k=0}(u^n_j)^{m-1-k}(U^n_j)^{k}e^n_j\bar{U}^n_j\right)_{\hat{x}}\right]\bar{e}^{n}_{j}\right|\nonumber\\
&=\left|-\frac{2bh}{m+2}\sum^{M-1}_{j=1}\left[(u^n_j)^{m}\bar{e}^n_j+\sum^{m-1}_{k=0}(u^n_j)^{m-1-k}(U^n_j)^{k}\bar{U}^n_je^n_j\right](\bar{e}^{n}_{j})_{\hat{x}}\right|\nonumber\\
&\leq Ch\sum^{M-1}_{j=1}\left[|\bar{e}^n_j|+|e^n_j|\right]|(\bar{e}^{n}_{j})_{\hat{x}}|\nonumber\\
&\leq C(\parallel e^{n+1}_{\hat{x}}\parallel^2+\parallel e^{n-1}_{\hat{x}}\parallel^2+\parallel e^{n+1}\parallel^2+\parallel e^{n}\parallel^2+\parallel e^{n-1}\parallel^2),\nonumber\\
&\leq C(\parallel e^{n+1}_x\parallel^2+\parallel e^{n-1}_x\parallel^2+\parallel e^{n+1}\parallel^2+\parallel e^{n}\parallel^2+\parallel e^{n-1}\parallel^2),
\end{align}
where theorem~\ref{theorem2}, Cauchy-Schwarz inequality and inequality (\ref{normin}) are used again.

Furthermore, we have
 \begin{align}\label{estimate3}
\textcolor{blue}{(r^{n},2\bar{e}^{n})\leq 2\parallel r^n\parallel \parallel \bar{e}^{n}\parallel \leq \parallel r^n\parallel^2+\parallel \bar{e}^{n}\parallel^2\leq\parallel r^n\parallel^2+\frac{1}{2}\left(\parallel e^{n+1}\parallel^2+\parallel e^{n-1}\parallel^2\right),}
\end{align}
\textcolor{blue}{where Cauchy-Schwarz inequality are used.}


Substituting (\ref{estimate1})-(\ref{estimate3}) into (\ref{estimate0}), we get
 \begin{align}\label{estimate4}
\parallel e^{n+1}\parallel^2&-\parallel e^{n-1}\parallel^2+\alpha(\parallel e^{n+1}_{x}\parallel^2-\parallel e^{n-1}_{x}\parallel^2)+\lambda(\parallel e^{n+1}_{x\bar{x}}\parallel^2-\parallel e^{n-1}_{x\bar{x}}\parallel^2)\nonumber\\
&\leq C\tau(\parallel e^{n+1}\parallel^2+\parallel e^{n}\parallel^2+\parallel e^{n-1}\parallel^2+\parallel e^{n+1}_{x}\parallel^2+\parallel e^{n-1}_{x}\parallel^2)+2\tau\parallel r^n\parallel^2.
\end{align}

Since $\alpha, \lambda$ are positive constants, it is easy to check that
 \begin{align*}
\parallel e^{n+1}\parallel^2&-\parallel e^{n-1}\parallel^2+\alpha(\parallel e^{n+1}_{x}\parallel^2-\parallel e^{n-1}_{x}\parallel^2)+\lambda(\parallel e^{n+1}_{x\bar{x}}\parallel^2-\parallel e^{n-1}_{x\bar{x}}\parallel^2)\nonumber\\
&\leq C\tau(\parallel e^{n+1}_{x}\parallel^2+\parallel e^{n}_{x}\parallel^2+\parallel e^{n+1}\parallel^2+\parallel e^{n}\parallel^2+\parallel e^{n-1}\parallel^2)+2\tau\parallel r^n\parallel^2\nonumber\\
&\leq C'\tau(\parallel e^{n+1}\parallel^2+2\parallel e^{n}\parallel^2+\parallel e^{n-1}\parallel^2+\alpha\parallel e^{n+1}_{x}\parallel^2+2\alpha\parallel e^{n}_{x}\parallel^2)+2\tau\parallel r^n\parallel^2,
\end{align*}
\textcolor{blue}{where $C'=\max(\frac{C}{\alpha}, C)$ and $C$ is the positive constant in the above inequality. }

\textcolor{blue}{Replacing $C'$ in the above inequality by the general positive constant notation $C$, we have}
\begin{align}\label{estimate6}
\parallel e^{n+1}\parallel^2&-\parallel e^{n-1}\parallel^2+\alpha(\parallel e^{n+1}_{x}\parallel^2-\parallel e^{n-1}_{x}\parallel^2)+\lambda(\parallel e^{n+1}_{x\bar{x}}\parallel^2-\parallel e^{n-1}_{x\bar{x}}\parallel^2)\nonumber\\
&\leq C\tau(\parallel e^{n+1}\parallel^2+2\parallel e^{n}\parallel^2+\parallel e^{n-1}\parallel^2+\alpha\parallel e^{n+1}_{x}\parallel^2+2\alpha\parallel e^{n}_{x}\parallel^2)+2\tau\parallel r^n\parallel^2\nonumber\\
&\leq C\tau(\parallel e^{n+1}\parallel^2+2\parallel e^{n}\parallel^2+\parallel e^{n-1}\parallel^2+\alpha\parallel e^{n+1}_{x}\parallel^2+2\alpha\parallel e^{n}_{x}\parallel^2+\alpha\parallel e^{n-1}_{x}\parallel^2\nonumber\\
&+\lambda\parallel e^{n+1}_{x\bar{x}}\parallel^2+2\lambda\parallel e^{n}_{x\bar{x}}\parallel^2+\lambda\parallel e^{n-1}_{x\bar{x}}\parallel^2)+2\tau\parallel r^n\parallel^2.
\end{align}


Let
\begin{align*}
D^n=\parallel e^{n}\parallel^2+\alpha\parallel e^{n}_{x}\parallel^2+\lambda\parallel e^{n}_{x\bar{x}}\parallel^2+\parallel e^{n-1}\parallel^2+\alpha\parallel e^{n-1}_{x}\parallel^2+\lambda\parallel e^{n-1}_{x\bar{x}}\parallel^2,
\end{align*}
then (\ref{estimate6}) can be rewritten as follows:
\begin{align*}
(D^{n+1}-D^n)\leq C\tau(D^{n+1}+D^n)+2\tau\parallel r^n\parallel^2,
\end{align*}
which is equivalent to
\begin{align}\label{estimate8}
(1-C\tau)(D^{n+1}-D^n)\leq 2C\tau D^n+2\tau\parallel r^n\parallel^2.
\end{align}

If $\tau$, which is sufficiently small, satisfies $\tau<\frac{1}{3C}$ ($C$ is the positive constant in the inequality (\ref{estimate8})), then $1-C\tau>0$ and (\ref{estimate8}) gives
 \begin{align*}
D^{n+1}-D^n &\leq \frac{2C}{1-C\tau}\tau D^n+\frac{2\tau}{1-C\tau}\parallel r^n\parallel^2\\
& \leq 3C\tau D^n+3\tau\parallel r^n\parallel^2,\\
& \leq C''(\tau D^n+\tau\parallel r^n\parallel^2),
\end{align*}
\textcolor{blue}{where $C''=\max(3C,3)$ and we have used $\frac{2}{1-C\tau}<3$ since $\tau<\frac{1}{3C}$.}

\textcolor{blue}{Replacing $C''$ in the above inequality by the general positive constant notation $C$, we have }
\begin{align}\label{estimate9}
D^{n+1}-D^n \leq C\tau D^n+C\tau\parallel r^n\parallel^2.
\end{align}

Summing (\ref{estimate9}) from $1$ to $n-1$, we get
 \begin{align}\label{estimate10}
D^{n}\leq D^1+C\tau \sum^{n-1}_{l=1} D^l+C\tau \sum^{n-1}_{l=1}\parallel r^l\parallel^2,
\end{align}
where
 \begin{align}\label{estimate11}
\tau \sum^{n-1}_{l=1}\parallel r^l\parallel^2\leq n\tau \max_{1\leq l\leq n-1}\parallel r^l\parallel^2\leq T\cdot O(\tau+h^2)^2.
\end{align}

Since $e^{0}_i=0$ and the Crank-Nicolson scheme (\ref{scheme4}) is used to compute $U^1$, we have $D^1=O(\tau^2+h^2)$ followed by a simple analysis for the scheme (\ref{scheme4}). Therefore
 \begin{align}\label{estimate1010}
D^{n}\leq O(\tau^2+h^2)^2+C\tau \sum^{n-1}_{l=1} D^l.
\end{align}

Using lemma~\ref{Lemma4}, we obtain
 \begin{align}\label{estimate12}
D^{n}\leq O(\tau^2+h^2)^2.
\end{align}

Thus,
 \begin{align}\label{estimate13}
\parallel e^n\parallel\leq O(\tau^2+h^2).
\end{align}

By using lemma~\ref{Lemma3}, we have
 \begin{align}\label{estimate14}
\parallel e^n\parallel_{\infty}\leq O(\tau^2+h^2),
\end{align}
i.e.,
 \begin{align}
\parallel u^n-U^n\parallel_{\infty}\leq C(\tau^2+h^2).
\end{align}
This completes  the proof.\\

\begin{theorem}\label{theorem5}
Suppose $u_0\in C^{7}_{0}[x_l,x_r]$, then the solution $U^n$ of the finite difference scheme (\ref{scheme1})-(\ref{scheme3}) is unconditionally stable with the $\parallel\cdot\parallel_{\infty}$ norm.
\end{theorem}

\textcolor{blue}{The proof of this theorem is similar as the above theorem.}

\section{Numerical results}\label{section4}

In order to preform numerical accuracy test of the method proposed in this paper, one needs to know the exact solitary solutions beforehand. \textcolor{blue}{Exact solitary solutions of the equation (\ref{KS1}) were previously obtained by~\cite{Zuon} using the sech ansatz method. In this paper, we provide a different derivation through the sine-cosine method, which is listed  in the appendix A.}

{\bf Example 1.} We present the numerical results for the case $m=2, a=1, b=1, c=2, \alpha=1, \lambda=1, \nu=1 $. From \textcolor{blue} {\cite{Zuon} and appendix A, one can obtain  the exact solitary solution as follows:}
 \begin{align}
u(x,t)=\frac{3}{4}\frac{\sqrt{370}-5\sqrt{10}}{\sqrt{5\sqrt{37}-29}}\textrm{sech}^2\left(\frac{\sqrt{\sqrt{37}-5}}{4}\left(x-\frac{33-5\sqrt{37}}{5\sqrt{37}-29}t\right)\right),
\end{align}
and the initial condition is set as
 \begin{align}\label{initial}
u(x,0)=\frac{3}{4}\frac{\sqrt{370}-5\sqrt{10}}{\sqrt{5\sqrt{37}-29}}\textrm{sech}^2\left(\frac{\sqrt{\sqrt{37}-5}}{4}x\right).
\end{align}

We first carry out the numerical convergence studies.  For the spatial convergence, we set $\tau=0.005$ as the fixed time step and use 4 different spatial meshes: $h=0.8, 0.4, 0.2, 0.1$, where $\tau$ is sufficient small such that the temporal error is negligible comparing to  the spatial error (Here the time step is $\tau=0.005$ while the smallest spatial size is $h=0.1$, thus $h^2>>\tau^2$, therefore, the dominant errors are the spatial errors). The final time $T$ is set to be $10$, and $x_l=-40, x_r=200$. Table \ref{table1} gives the errors between numerical solutions and exact solutions. We can see that the error decreases when the spatial mesh is refined and the convergence rate is two. Thus, the method is second-order convergent in space, which is consistent with theoretical results in the above section. For the temporal convergence, we set $h=0.005$ as the fixed spatial mesh and use 4 different temporal meshes: $\tau=0.8, 0.4, 0.2, 0.1$,  where $h$ is sufficient small such that the spatial error is negligible comparing to the temporal  error (Here the smallest time step is $\tau=0.1$ while the spatial size is $h=0.005$, thus $h^2<<\tau^2$, therefore, the dominant errors are the temporal errors). The final time $T$ is set to be $10$, and $x_l=-40, x_r=200$. Table \ref{table2} gives the errors between numerical solutions and exact solutions. Again,  we can see that the error decreases when the temporal mesh is refined, and the convergence rate is also two. Thus, the method is second-order convergent in time, which is again consistent with theoretical results in the above section.

In order to show that the numerical scheme has the energy conservative property (\ref{disenergy}),
we carry out another computation, where $T=100, x_{l}=-40, x_{r}=240, h=0.1, \tau=0.1$ are used. Table \ref{table3}
gives the quantities of $E^{n}$ at several time stages.  We can see that $E^{n}$ is conserved exactly (up to 8 decimals) during the time evolution of the solitary wave. 

{\bf Example 2.} We present the numerical results for the case $m=4, a=1, b=1, c=2, \alpha=1, \lambda=1, \nu=1 $.  From \textcolor{blue} {\cite{Zuon} and appendix A, one can obtain  the exact solitary solution as follows:}
 \begin{align}
u(x,t)=\left(\frac{40(\sqrt{127}-10)^2}{3(10\sqrt{127}-109)}\right)^{\frac{1}{4}}\textrm{sech}\left(\frac{\sqrt{\sqrt{127}-10}}{3}\left(x-\frac{118-10\sqrt{127}}{10\sqrt{127}-109}t\right)\right),
\end{align}
and the initial condition is set as
 \begin{align}\label{initial2}
u(x,0)=\left(\frac{40(\sqrt{127}-10)^2}{3(10\sqrt{127}-109)}\right)^{\frac{1}{4}}\textrm{sech}\left(\frac{\sqrt{\sqrt{127}-10}}{3}x\right).
\end{align}

Again, we carry out the spatial and temporal convergence. \textcolor{blue}{Table \ref{table4} and \ref{table5}} give the errors between numerical solutions and exact solutions for spatial and temporal convergence, respectively.  Once again, we can see that the method is second-order convergent both in time  and space variables.

Additionally, table \ref{table6}
provides the quantities of $E^{n}$  with $T=100, x_{l}=-40, x_{r}=240, h=0.1, \tau=0.1$. Once again, we can see that $E^{n}$ is conserved and the  method can be well used to study the solitary wave at long time.

\section{Conclusions }
In  this paper, a three-level linearly implicit  finite difference method for the initial value problem of the generalized Rosenau-Kawahara-RLW equation is developed. The fundamental energy conservative property is preserved by the current numerical scheme. The existence and uniqueness of the numerical solution are also proved.  The method is shown to be second-order  convergent both in time and space variables,  and it is  unconditionally stable. Moreover, exact solitary solutions are derived through sine-cosine method which is used for the numerical  tests.  Numerical results confirm well with the theoretical results.

\section*{Acknowledgement}
Dongdong He was supported by the Program for Young Excellent Talents at Tongji University (No. 2013KJ012), the Natural Science Foundation of China (No. 11402174) and the Scientific Research Foundation for the Returned Overseas Chinese Scholars, State Education Ministry. Kejia Pan was supported by the Natural Science Foundation of China (Grant Nos. 41474103 and 41204082), the Natural Science Foundation of Hunan Province of China (Grant No. 2015JJ3148) and the Mathematics and Interdisciplinary Sciences Project of Central
South University.

\clearpage
\begin{table}[b]
\caption{Spatial mesh refinement analysis with $\tau=0.005, T=10$ for example 1.}\label{table1}
\begin{center}
\begin{tabular}{|c|c|c|c|c|}
\hline
$h$&$\parallel e\parallel$&rate&$\parallel e\parallel_{\infty}$&rate\\
\hline $0.8  $&  2.66e-1 &    -&1.032e-1&   -\\
\hline $0.4  $&6.650-2 &     2.002    & 2.570-2  &      2.006  \\
\hline $0.2  $&1.666-2  & 1.997 &  6.460-3 &    1.992 \\
\hline $0.1 $& 4.209-3 &   1.985 & 1.631-3&   1.986
\\
\hline
\end{tabular}
\end{center}
\end{table}

\clearpage
\begin{table}[b]
\caption{Temporal mesh refinement analysis with $h=0.005, T=10$ for example 1.}\label{table2}
\begin{center}
\begin{tabular}{|c|c|c|c|c|}
\hline
$\tau$&$\parallel e\parallel$&rate&$\parallel e\parallel_{\infty}$&rate\\
\hline 0.8& 1.523 &    -&6.060e-1&   -\\
\hline 0.4 &3.812e-1 &     1.999  &1.540e-1&   1.976\\
\hline 0.2 &9.711e-2 &    1.973  &3.910e-2& 1.977 \\
\hline 0.1 &2.442e-2 &    1.992  &9.820e-3&  1.993\\
\hline
\end{tabular}
\end{center}
\end{table}
\clearpage

\begin{table}[b]
\caption{Invariant of $E^{n}$ for example 1.}\label{table3}
\begin{center}
\begin{tabular}{|c|c|}
\hline
$t$&$E^{n}$\\
\hline 0.05&   25.451405792697514\\
\hline 19.95&  25.451405792693116  \\
\hline 39.95&  25.451405792447929\\
\hline 59.95&  25.451405792214793\\
\hline 79.95&  25.451405791920855 \\
\hline 99.95&  25.451405792207414\\
\hline
\end{tabular}
\end{center}
\end{table}
\clearpage

\begin{table}[b]
\caption{Spatial mesh refinement analysis with $\tau=0.005, T=10$ for example 2.}\label{table4}
\begin{center}
\begin{tabular}{|c|c|c|c|c|}
\hline
$h$&$\parallel e\parallel$&rate&$\parallel e\parallel_{\infty}$&rate\\
\hline $0.8  $& 1.543e-1 &    -&5.839e-2&   -\\
\hline $0.4  $&3.790-2 &     2.026   & 1.446e-2  &      2.013 \\
\hline $0.2  $&9.440-3  &   2.005 &  3.599e-3 &  2.007\\
\hline $0.1 $& 2.366-3 &  1.996 &9.011e-4&  1.999
\\
\hline
\end{tabular}
\end{center}
\end{table}

\clearpage
\begin{table}[b]
\caption{Temporal mesh refinement analysis with $h=0.005, T=10$ for example 2.}\label{table5}
\begin{center}
\begin{tabular}{|c|c|c|c|c|}
\hline
$\tau$&$\parallel e\parallel$&rate&$\parallel e\parallel_{\infty}$&rate\\
\hline 0.8& 4.582-1 &    -&2.029e-1&   -\\
\hline 0.4 &8.633e-2 &      2.408  &3.843e-2&     2.401 \\
\hline 0.2 &2.124e-2 &     2.023    &9.447e-3&    2.024 \\
\hline 0.1 &5.263e-3 &     2.012  &2.330e-3&   2.019\\
\hline
\end{tabular}
\end{center}
\end{table}
\clearpage

\begin{table}[b]
\caption{Invariant of $E^{n}$ for example 2.}\label{table6}
\begin{center}
\begin{tabular}{|c|c|}
\hline
$t$&$E^{n}$\\
\hline 0.05&   13.565665615099391\\
\hline 19.95&  13.565665614771643\\
\hline 39.95&  13.565665614965912\\
\hline 59.95&  13.565665614937172\\
\hline 79.95&  13.565665614960499\\
\hline 99.95& 13.565665614998375\\
\hline
\end{tabular}
\end{center}
\end{table}
\clearpage

\appendix
\section{Exact solitary solutions}

The sine-cosine method uses the sine or cosine function as the wave function to seek the traveling wave solution of a time dependent partial differential equation, which has the advantage of reducing the nonlinear problem to a
system of algebraic equations that can be easily solved by using a symbolic computation system such as Mathematica or Maple~\cite{Kudryashov,Zuo1,Mdallal,Wongsaijai}.

For equation (\ref{KS1}), one can obtain the exact solitary solutions by using the sine-cosine method. Firstly, we seek the following travelling wave solution~\cite{Kudryashov,Zuo1,Mdallal,Wongsaijai}:
 \begin{align}\label{traveling}
u(x,t)=\hat{u}(\xi),\quad  \xi=x-vt,
\end{align}
where $v$ is referred as the wave velocity which is a constant to be determined later.

Under the transformation of (\ref{traveling}) and integrating once,  (\ref{KS1}) can be reduced into:
 \begin{align}\label{utranving}
(a-v)\hat{u}+\frac{b}{m+1}\hat{u}^{m+1}+(v\alpha+c)\hat{u}_{\xi\xi}-(\lambda v+\nu) \hat{u}_{\xi\xi\xi\xi}=0.
\end{align}
Secondly, by using the sine-cosine method, the solutions of the above reduced ODE equation can be expressed in the following form~\cite{Zuo1,Mdallal,Wongsaijai}:
\begin{equation}\label{form1}
\hat{u}(\xi) = \left\{ \begin{array}{ll}
A\cos^{\eta}(B\xi),\ \ \ \ \ & \textrm{if $|\xi|<\frac{\pi}{2B}$},\\
0,\ \ \ \ & \textrm{otherwise},
\end{array} \right.
\end{equation}
or in the form
\begin{equation}
\hat{u}(\xi) = \left\{ \begin{array}{ll}
A\sin^{\eta}(B\xi), \ \ \ \ \ & \textrm{if $|\xi|<\frac{\pi}{2B}$},\\
0,\ \ \ \ & \textrm{otherwise},
\end{array} \right.
\end{equation}
where $A, B, \eta$ are parameters to be determined.

Using (\ref{form1}), one have
 \begin{align}\label{u2}
\hat{u}_{\xi\xi}=AB^2\eta(\eta-1)\cos^{\eta-2}(B\xi)-AB^2\eta^2\cos^{\eta}(B\xi),
\end{align}
and
 \begin{align}\label{u4}
\hat{u}_{\xi\xi\xi\xi}=AB^4\eta(\eta-1)(\eta-2)(\eta-3)\cos^{\eta-4}(B\xi)-2AB^4\eta(\eta-1)(\eta^2-2\eta+2)\cos^{\eta-2}(B\xi)+AB^4\eta^4\cos^{\eta}(B\xi).
\end{align}
Substituting (\ref{form1}), (\ref{u2}) and (\ref{u4}) into (\ref{utranving}), one obtain
\begin{align}
&\frac{bA^{m+1}}{m+1}\cos^{\eta (m+1)}(B\xi)-AB^4\eta(\eta-1)(\eta-2)(\eta-3)(\lambda v+\nu)\cos^{\eta-4}(B\xi)+\nonumber\\
&(AB^2(v\alpha+c)\eta(\eta-1)+2AB^4\eta(\eta-1)(\eta^2-2\eta+2)(\lambda v+\nu))\cos^{\eta-2}(B\xi)+\nonumber\\
&\left((a-v)A-AB^2(v\alpha+c)\eta^2-AB^4\eta^4(\lambda v+\nu)\right)\cos^{\eta}(B\xi)=0.
\end{align}
Balancing $\cos^{\eta (m+1)}(B\xi)$ and $\cos^{\eta-4}(B\xi)$, and setting each coefficients
of $\cos^{j}(B\xi)$ ($j=\eta, \eta-2, \eta-4$)  to be zero, one can obtain a set of equations for $\eta, v, A, B$ as follows:
 \begin{align}
&\eta (m+1)=\eta-4,\label{eta}\\
&\frac{bA^{m+1}}{m+1}-AB^4\eta(\eta-1)(\eta-2)(\eta-3)(\lambda v+\nu)=0,\label{eq2}\\
&AB^2(v\alpha+c)\eta(\eta-1)+2AB^4\eta(\eta-1)(\eta^2-2\eta+2))(\lambda v+\nu)=0,\label{eq3}\\
&(a-v)A-AB^2(v\alpha+c)\eta^2-AB^4\eta^4(\lambda v+\nu)=0.\label{eq4}
\end{align}

From (\ref{eta}), one can  find that
 \begin{align}
\eta=-\frac{4}{m}.\label{eta0}
\end{align}

From (\ref{eq3}) and (\ref{eq4}), one obtain that
\begin{align}
B^2=\frac{(\lambda a+\nu)(\eta^2-2\eta+2)\pm\sqrt{(\lambda a+\nu)^2(\eta^2-2\eta+2)^2+(\lambda c-\nu\alpha)(\alpha a+c)\eta^2(\eta-2)^2}}{(\lambda c-\nu\alpha)\eta^2(\eta-2)^2}.\label{B1}
\end{align}

Once $B$ is obtained, one can get that
 \begin{align}
&v=-\frac{2\nu B^2(\eta^2-2\eta+2)+c}{2\lambda B^2(\eta^2-2\eta+2)+\alpha}, \label{v}\\
&A=\left[\frac{m+1}{b}B^4\eta(\eta-1)(\eta-2)(\eta-3)(\lambda v+\nu)\right]^{\frac{1}{m}}. \label{A}
\end{align}

The classification for the above solution is discussed as follows:

Case 1: $(\lambda a+\nu)^2(\eta^2-2\eta+2)^2+(\lambda c-\nu\alpha)(\alpha a+c)\eta^2(\eta-2)^2>0$.

(1):  $\frac{(\lambda a+\nu)(\eta^2-2\eta+2)+\sqrt{(\lambda a+\nu)^2(\eta^2-2\eta+2)^2+(\lambda c-\nu\alpha)(\alpha a+c)\eta^2(\eta-2)^2}}{(\lambda c-\nu\alpha)\eta^2(\eta-2)^2}>0$.
We obtain a periodic solution as follows:
 \begin{align}
u(x,t)=\hat{u}(\xi)=A\cos^{\eta}(B(x-vt)),\label{exact00}
\end{align}
 where $B$ is given by
  \begin{align}
B=\sqrt{\frac{(\lambda a+\nu)(\eta^2-2\eta+2)+\sqrt{(\lambda a+\nu)^2(\eta^2-2\eta+2)^2+(\lambda c-\nu\alpha)(\alpha a+c)\eta^2(\eta-2)^2}}{(\lambda c-\nu\alpha)\eta^2(\eta-2)^2}},
\end{align}
and $\eta, v, A$ are given by (\ref{eta0}), (\ref{v}) and (\ref{A}), respectively.

(2): $\frac{(\lambda a+\nu)(\eta^2-2\eta+2)+\sqrt{(\lambda a+\nu)^2(\eta^2-2\eta+2)^2+(\lambda c-\nu\alpha)(\alpha a+c)\eta^2(\eta-2)^2}}{(\lambda c-\nu\alpha)\eta^2(\eta-2)^2}<0$. We obtain a solitary solution as follows:
\begin{align}
u(x,t)=\hat{u}(\xi)=A\textrm{sech}^{\eta}(B_0(x-vt)),\label{exact0}
\end{align}
where $B_0$ is given by
\begin{align}
B_0=\sqrt{\frac{-(\lambda a+\nu)(\eta^2-2\eta+2)-\sqrt{(\lambda a+\nu)^2(\eta^2-2\eta+2)^2+(\lambda c-\nu\alpha)(\alpha a+c)\eta^2(\eta-2)^2}}{(\lambda c-\nu\alpha)\eta^2(\eta-2)^2}},
\end{align}
and  $\eta, v, A$ are given by (\ref{eta0}), (\ref{v}) and (\ref{A}), respectively.

(3):  $\frac{(\lambda a+\nu)(\eta^2-2\eta+2)-\sqrt{(\lambda a+\nu)^2(\eta^2-2\eta+2)^2+(\lambda c-\nu\alpha)(\alpha a+c)\eta^2(\eta-2)^2}}{(\lambda c-\nu\alpha)\eta^2(\eta-2)^2}>0$. Similar periodic solution as (\ref{exact00}) can be obtained.

(4): $\frac{(\lambda a+\nu)(\eta^2-2\eta+2)-\sqrt{(\lambda a+\nu)^2(\eta^2-2\eta+2)^2+(\lambda c-\nu\alpha)(\alpha a+c)\eta^2(\eta-2)^2}}{(\lambda c-\nu\alpha)\eta^2(\eta-2)^2}<0$. Similar solitary solution as (\ref{exact0}) can be obtained.

Case 2: $(\lambda a+\nu)^2(\eta^2-2\eta+2)^2+(\lambda c-\nu\alpha)(\alpha a+c)\eta^2(\eta-2)^2<0$.

In this case, the resulting $A, B, v$ are complex numbers with nonzero real parts.  After substituting the expressions of $A, B, v, \eta$  into equation (\ref{form1}), $u$ needs to take the real part or the imaginary part of the resulting equation. The detail is omitted here.

Similarly, one can use the sine method  to get another set of solutions.\\

\end{document}